\documentclass[a4paper,11pt]{article}
\usepackage{amssymb}

\usepackage{amsthm}
\usepackage{pstricks}
\usepackage{pstcol}
\usepackage{amsmath}
\usepackage{epsfig}
\usepackage{wrapfig}

\theoremstyle{definition}
\theoremstyle{remark}

\begin{document}

\date{}
\title{\textbf{On a Watson-like Uniqueness Theorem and Gevrey Expansions.}}
\author{\textbf{\ \ D.W.H. Gillam~${^{\text{}}}$, {V. Gurarii}~}${^{\mathbf{}}}
\footnote{This work was initiated in the School of Mathematical
Sciences of Swinburne University of Technology. One of us (VG) is
grateful to the Weizmann Institute of Science, Israel, for support
under the Rosi and Max Varon Visiting Professorship program during
the last stage of this work.
\newline
\textsl{Math. classification 2000:34M25, 34M30, 34M37, 34M40}}$}
\maketitle

\begin{abstract}
We present a maximal class of analytic functions, elements of which are 
in one-to-one correspondence with their asymptotic expansions.
In recent decades it has been realized (B. Malgrange, J. Ecalle,
J.-P. Ramis, Y. Sibuya et al.), that the formal power series
solutions of a wide range of systems of ordinary (even non-linear)
analytic differential equations are in fact the Gevrey expansions
for the regular solutions. Watson's uniqueness theorem belongs to
the foundations of this new theory. This paper contains a
discussion of an extension of Watson's uniqueness theorem for
classes of functions which admit a Gevrey expansion in angular
regions of the complex plane with opening less than or equal to
\(\frac \pi k, \) where \(k\) is the order of the Gevrey
expansion. We present conditions which ensure uniqueness and which
suggest an extension of Watson's representation theorem. These
results may be applied for solutions of certain classes of
differential equations to obtain the best accuracy estimate for
the deviation of a solution from a finite sum of the corresponding
Gevrey expansion.

\end{abstract}

%%%%%%%%%%%%%%%%%%%%%%%%%%%%%%%%%%%%%%%%%%%%%%%%%%%%%%%%

\begin{itemize}
\item  \(\mathbb C\) stands for the complex plane;

\item  \(S(\alpha ,\beta )\) stands for the sector
\[
S(\alpha ,\beta )=\left\{ z\in \mathbb{C}:\,0<\left| z\right| <\infty
,\alpha <\arg z<\beta \right\} \label{Sect}
\]
where the number \(\beta -\alpha \) is said to be the
\textsl{opening} of the sector \(S(\alpha ,\beta) \).
\end{itemize}

\textbf{Introduction.} In 1912 G. Watson published the following
result, \cite{Wat1}, section 8, Theorem V.

\vspace*{2.5ex}

\textbf{Watson's uniqueness theorem.}\textsl{\ Let
}\(\{p_0,p_1,\ldots \}\) \textsl{\ be a given sequence of complex
numbers, and }\(P(z)\)\textsl{\ be a function satisfying the
conditions:}

\begin{enumerate}
\item[\textup{(i)}]  \({P}\left( z\right) \)\textsl{\ is analytic and
single-valued in the sector }\(S\left( \alpha ,\beta \right)
\)\textsl{;}

\item[\textup{(ii)}]  \({P}\left( z\right) \)\textsl{\ admits the following
series of estimates for all \(z\in S\left(
\alpha ,\beta \right)\) except \(|z|<\sigma\) }
\begin{equation}
\left| P\left( z\right) -\sum_{k=0}^{n-1}\frac{p_k}{z^{k+1}}\right| \leq M%
\frac{n!}{a^n\left| z\right| ^{n+1}},\ \ n=0,\,1,\,\ldots ,  \label{E_S}
\end{equation}
\textsl{where positive constants \(M,\,\,a\) and \(\,\sigma\,\)%
 do not depend on \(z\) or }\(n.\)
\end{enumerate}

\textsl{If the opening of the sector }\(S\left( \alpha ,\beta
\right) \) \textsl{\ satisfies the condition }\(\beta -\alpha >\pi
\)\textsl{, then the function }\(P(z)\)\textsl{\ is uniquely
determined by conditions (i) and (ii): two functions
}\(P_1(z)\)\textsl{\ and }\(P_2(z)\)\textsl{\ satisfying
conditions (i) and (ii) with the same sequence }\(\{p_0,p_1,\ldots \}\)\textsl{%
\ must coincide on }\(S\left( \alpha ,\beta \right) \)\textsl{.}

\vspace{2.5ex}

Set
\begin{equation} \hat{P}(z)=\sum_{k=0}^{\infty
}\frac{p_{k}}{z^{k+1}} .  \label{fP}
\end{equation}

The formal power series \(\hat{P}(z)\) is known as \textsl{the
Gevrey expansion of order \(1\)} for a function \({P}(z)\)
satisfying conditions (i) and (ii).

In recent times it has been realized (B. Malgrange, J. Ecalle,
J.-P. Ramis, Y. Sibuya et al.), that the formal power series
solutions of a wide range of systems of ordinary (even non-linear)
analytic differential equations are in fact the Gevrey expansions
for the regular solutions, see, for example, \cite{Ramis1},
\cite{RamisSi} and \cite{RamisSi1}. Watson's uniqueness theorem
belongs to the foundations of this new theory.

\vspace{2.5ex}

To prove Watson's uniqueness theorem we introduce the function
\[
P(z)=P_1(z)-P_2(z)
\]
which satisfies the inequalities
\begin{equation}
\left| P\left( z\right) \right| \leq 2M\frac{n!}{a^n\left| z\right| ^{n+1}}%
,\ \ n=0,\,1,\,\ldots ,\ \ z\in S\left( \alpha ,\beta \right) .  \label{EP_S}
\end{equation}
Minimizing the right-hand side of (\ref{EP_S}), for fixed \(z,\,
|z|>\frac{1}{a} \), with respect to \(n\) yields the inequality
\begin{equation}
\left| P\left( z\right) \right| \leq M_a e^{-a|z|} \label{P_S}
\end{equation}
where \(M_a=4M\sqrt{2\pi}a.\)

Thus, Watson's uniqueness theorem can be derived from the
following fact, a proof of which may be found in \cite{Hardy},
section 8.11.

\vspace*{1.5ex}

\textbf{Lemma 1}. \textsl{Let }\(P\left( z\right) \)\textsl{\ be
an analytic function in the sector\newline }\({S_\varepsilon
}=S\left( -\frac \pi 2\left( 1+\varepsilon \right) ,\frac \pi
2\left( 1+\varepsilon \right) \right) ,\)\textsl{\
}\(0<\varepsilon <1,\) \textsl{\ satisfying the following estimate
}
\begin{equation}
\left| P\left( z\right) \right| \leq Me^{-a\left| z\right| }\,\,\text{for all%
}\,\,\,z\in {S_\varepsilon },  \label{E_PM}
\end{equation}
\textsl{where }\(a\)\textsl{\ is a positive constant. Then
}\(P\left( z\right) \equiv 0.\)

\vspace*{1.5ex}

We note that, since the estimate (\ref{P_S}) holds for \(z\in
S\left( \alpha ,\beta \right) \backslash \left\{ z:\left| z\right|
<\frac 1a\right\}\), it holds also in the shifted sector \(S\left(
\alpha ,\beta \right)+de^{i\frac{\alpha+\beta}{2}}\) where
\(d=\frac{1}{a\cos{\frac {\pi\varepsilon}{2} }}.\) Thus, strictly
speaking, lemma 1 should be applied rather to the shifted function
\(P\left( \left(z+d\right)e^{i\frac{\alpha +\beta }2}\right).\)

\vspace{1.5ex}

If the opening of the sector \(S\left( \alpha ,\beta \right) \)
satisfies the condition \(\beta -\alpha <\pi \), then the above
uniqueness result is indisputably false. Indeed, given
\(\delta ,\,0<\delta <\frac \pi 2 ,\) every function \(P\left(
z\right) \) of the form
\begin{equation}
P\left( z\right) =\frac{\varphi \left( z\right) }ze^{-z},  \label{fi}
\end{equation}
where \(\varphi \left( z\right) \) is analytic and bounded in the
sector \(S\left( -\frac \pi 2+\delta ,\frac \pi 2-\delta \right)
\), satisfies the conditions (\ref{E_S}) with
\[
M={\sup_{{S\left( -\frac \pi 2+\delta ,\frac \pi 2-\delta \right) }}}\left|
\varphi \left( z\right) \right| ,
\]
\(a=\sin \delta \) and \(p_k=0,k=0,1,\ldots \). This can be easily
derived using the elementary inequalities
\[
e^{-\left| z\right| }<\frac{n!}{\left| z\right| ^n},\,\,n=0,\,1,\,\ldots \text{.}
\]
Thus, in this case the set of functions which satisfy conditions
(i) and (ii) of Watson's theorem with \(\{p_0,p_1,\ldots
\}=\{0,0,\ldots \}\) is rather large.

The two opposed cases, \(\beta -\alpha >\pi \) and \(\beta -\alpha
<\pi ,\)
discussed above lead us to regard the case \(\beta -\alpha =\pi \) as the%
\textsl{\ critical case} and the value \(\pi \) as the
\textsl{critical value} of the opening \(\beta -\alpha \). If the
opening of \(S\left( \alpha ,\beta \right) \) takes this critical
value the situation is much more delicate.

However, the statement of Lemma 1 can be extended to this critical
case.

\vspace{1.0ex}

\textbf{Lemma 2.}\textsl{\ Let }$P\left( z\right) $\textsl{\ be an analytic
function in the sector }$S\left( -\frac \pi 2,\frac \pi 2\right) $\textsl{\
satisfying the following estimate }
\begin{equation}
\left| P\left( z\right) \right| \leq Me^{-a\left| z\right| }\,\,\text{for all%
}\,\,z\in S\left( -\frac \pi 2,\frac \pi 2\right) ,  \label{PM}
\end{equation}
\textsl{where }\(a\)\textsl{\ is a positive constant. Then
}\(P\left( z\right) \equiv 0.\)

\vspace{1.0ex}

\textsf{Proof of Lemma 2}. We demonstrate three different proofs
of Lemma 2.

(1) Although Hardy's argument in \cite{Hardy}, section 8.11 
for the proof of lemma 1 is not
applicable in this case it is still possible to apply the
Phragm\'{e}n-Lindel\"{o}f theorem to prove the result, and the
proof can be made even simpler. Given \(k,\,k>0,\) consider the
function \(P_k(z)=P(z)e^{kz}\) inside the sector
\(S(-\arctan{\frac k a},\arctan{\frac k a})\) where \(P\left(
z\right) \) is a function satisfying the conditions of lemma 2. On
the boundary of this sector we have \(a|\Im z|-k\Re z=0\) and so
\(|P_k(z)| \leq M \). The standard Phragm\'{e}n-Lindel\"{o}f
theorem yields the same estimate inside the sector. Thus, for
\(z>0\) we have \(|P(z)| \leq M e^{-kz}\). As  \(k\) is
independent of \(z\) we can let \(k\) tend to \(+ \infty \) to
obtain \(P\left( z\right) \equiv 0.\) \hfill \(\blacktriangle \)

(2) Alternatively, the condition (\ref{PM}) allows one to use the
following well-known theorem of complex analysis.

\vspace{1.0ex}

\textbf{Theorem}.\,\footnote{%
This theorem follows from the so-called Jensen inequality on the
summability of the logarithmic integral of a bounded analytic
function. The integral on the left-hand side of (\ref{log})
generated an extensive literature, culminating in the two-volume
monograph by Paul Koosis, The logarithmic integral, I and II,
Cambridge Study in Advance Mathematics, Volumes 12 and 21,
Cambridge University Press, (1992), (1998).} \textsl{\ If
}\(P\left( z\right) \)\textsl{\ is analytic and bounded in the
sector }\(S\left( -\frac \pi 2,\frac \pi 2\right) \)\textsl{\ and
for some }\(c>0\)\textsl{\ }
\begin{equation}
\int_{c-i\infty }^{c+i\infty }\frac{\log \left| P\left( z\right) \right| }{%
1+\left| z\right| ^2}d\left| z\right| =-\infty ,  \label{log}
\end{equation}
\textsl{then }\(P\left( z\right) \equiv 0.\)

\vspace{1.0ex}

Indeed, if \(P\left( z\right) \) satisfies (\ref{PM}) then
\(P\left( z\right) \) satisfies (\ref{log}).\hfill
\(\blacktriangle \)

(3) It is also possible to use other arguments to prove Lemma 2;
the idea of these arguments will be employed and extended later to
prove the main result of our paper. Let \(P\left( z\right) \) be a
function satisfying the conditions of lemma 2. Introduce the
function

\begin{equation}
F\left( t\right) =\frac 1{2\pi i}\int_{-i\infty }^{i\infty
}P\left( z\right) e^{tz}dz. \label{Bo}
\end{equation}

Condition (\ref{PM}) ensures the absolute convergence of the
integral in (\ref{Bo}) for every complex \(t,\,-a<\Im t<a\). Thus
\(F\left( t\right)\) is an analytic function in the strip
\(\{t:\left| \Im t\right| <a\}\) of \(\mathbb{C}_t\). On the other
hand, for \(t \in (- \infty , a) \) the function
\(P_t(z)=P(z)e^{tz}\) is analytic and bounded, \(|P_t(z)| \leq M,
\) in the half-plane \(\{z:\Re z>0\}\) of \(\mathbb{C}_z\).
Moreover, \(P(z)e^{tz}\rightarrow 0\) as \(|z|\rightarrow
\infty.\) Using the Cauchy theorem this yields \( F\left( t\right)
=0\) for \(t<a\). From the uniqueness theorem for analytic
functions it follows that \(F\left( t\right) \equiv 0\). Thus,
\(P\left( z\right) \equiv 0.\)\hfill \(\blacktriangle \)

\vspace{1.0ex}

Lemma 2 shows that the conclusion of Watson's uniqueness theorem,
as   stated above, still remains valid for the critical sector,
that is the sector \(S(\alpha,\beta)\) with the critical values of
its opening satisfying \(\beta-\alpha=\pi\).

However, and this will be explained below, a certain weakening of
condition (ii) of Watson's uniqueness theorem leads to a loss of
the uniqueness property.

\vspace{1.0ex}

\textbf{Comment 1}. It is worthwhile noting that Watson's
uniqueness theorem can be transplanted to any sector \(S(\alpha
,\beta )\) of the Riemann surface of \(\log z\) with opening
greater than or equal to \(\frac \pi k\), where \(k\) is a
positive number. The following extension of Watson's uniqueness
theorem is valid.

\vspace{1.5ex}

\textbf{Watson's uniqueness theorem for order
\boldmath{\(k\)}.}\textsl{\ Let \(k\) be a positive number. Let
\(\{p_0,p_1,\ldots \}\) be a given sequence of complex numbers,
and \(P(z)\) a function satisfying the conditions:}

\begin{enumerate}
\item[\textup{(i)}]  \({P}\left( z\right) \)\textsl{\ is analytic and
single-valued in the sector }\(S\left( \alpha ,\beta \right)
\)\textsl{;}

\item[\textup{(ii)}]  \({P}\left( z\right) \)\textsl{\ admits the following
series of estimates for all \(z\in S\left(
\alpha ,\beta \right)\) except \(|z|<\sigma\)}
\begin{equation}
\left| P\left( z\right) -\sum_{j=0}^{n-1}\frac{p_j}{z^{j+1}}\right| \leq M%
\frac{\left( n!\right) ^{\frac 1k}}{\left( ka\left| z\right| \right) ^{n+1}}%
,\ \ n=0,\,1,\,\ldots ,\,\,\text{for all},  \label{k-estimates}
\end{equation}
\textsl{where positive constants \(M,\,\,a\) and \(\,\sigma\,\)%
 do not depend on \(z\) or }\(n.\)
\end{enumerate}

\textsl{If the opening of the sector }\(S\left( \alpha ,\beta
\right) \) \textsl{\ satisfies the condition }\(\beta -\alpha \geq
\frac \pi k\)\textsl{, then the function }\(P(z)\)\textsl{\ is
uniquely determined by the conditions (i) and (ii).}%
\hspace*{0.01ex}\hfill \(\blacktriangle\)

\vspace{1.0ex}

This theorem follows immediately from the corresponding extension
of Lemma 2 using the map \(\zeta =z^k.\) Comparing again the two
opposed cases \(\beta -\alpha >\frac \pi k\) and \(\beta -\alpha
<\frac \pi k\) shows that for this situation a sector \(S\left(
\alpha ,\beta \right) \) of the Riemann surface of \(\log z\) will
be a critical sector if its opening \(\beta -\alpha \) is equal to
\(\frac \pi k.\)

Expansions satisfying conditions similar to (\ref{k-estimates})
are known as Gevrey expansions of order $k$, (see \cite{Gevrey},
\cite{Ramis} and \cite{Ramis1}.)

\vspace{1.0ex}

Watson's uniqueness theorem was also extended by T. Carleman, see
\cite{Ca}. Carleman replaced \(n!\) by a sequence of positive
numbers \(m_n\) and, under certain regularity conditions on the
growth of \(m_n\) as \(n\rightarrow \infty \), he gave necessary
and sufficient conditions for uniqueness.

\vspace{1.5ex}

\textbf{The main theorem.} We develop Watson's result in a
different direction, keeping in mind the possible application of
this development to the analytic theory of differential equations
in the complex plane. Given \(\delta , \, 0<\delta <\frac \pi 2,\)
we consider a function \(P(z)\) which satisfies almost all the
conditions of Watson's uniqueness theorem, the only difference
being that in  (\ref{E_S}) the condition \(\, z\in S\left( \alpha
,\beta \right) \) is replaced by \(z\in S\left( -\frac \pi
2+\delta ,\frac \pi 2-\delta \right)\) and the positive constant
\(M\) in (\ref{E_S}) is replaced by a positive function
\(M(\delta) \) which may depend on \(\delta \). We plan to
consider later the case in which the constant \(a\) in (\ref{E_S})
is  replaced by a positive function \(a(\delta) \) which may also
depend on \(\delta \).

Consider a positive function \(M\left( \delta \right),\) defined on
the interval \(0<\delta <\frac \pi 2,\) where, for simplicity, we
assume that \(\log M\left( \delta \right)>1.\) Our aim is to
clarify under what conditions
 on \(M(\delta) \) the set of
functions satisfying (\ref{E_S}) in every sub-sector \(S\left(
-\frac \pi 2+\delta ,\frac \pi 2-\delta \right)\) of \( S\left(
-\frac \pi 2,\frac \pi 2\right) \), with \(M=M(\delta) \),
possesses a uniqueness property as stated in Watson's uniqueness
theorem.

The following extension of Watson's uniqueness theorem for the
critical sector \(S\left( -\frac \pi 2,\frac \pi 2\right) \) is
valid.

\vspace{1.0ex}

\textbf{Theorem 1.} \textsl{Let \(\{p_0,p_1,\ldots \}\)  be a
given sequence of complex numbers, and \(P(z)\) be a function
satisfying the conditions:}

\begin{enumerate}
\item[\textup{(i)}]  \textsl{\({P}\left( z\right) \) is analytic and
single-valued in the sector \( S\left( -\frac \pi 2,\frac \pi
2\right) \);}

\item[\textup{(ii)}]  \textsl{For every \(\delta ,0<\delta <\frac \pi 2,\)
\({P}\left( z\right) \) admits the following series of estimates}
\begin{align}
|P(z)&
-\sum_{k=0}^{n-1}\frac{p_k}{z^{k+1}}| \leq K_PM(\delta)%
\frac{n!}{a^n| z| ^{n+1}},\ \ n=0,\,1,\,\ldots ,\label{Est_KP}\\[1.0ex]
{}& \text{for all}\,\,\,\,z\in S( -\frac \pi 2+\delta ,\frac \pi
2-\delta ) ,\,\,|z|>\sigma,\notag
\end{align}

\textsl{where \(a,\,\,\sigma\) and \(K_P\)\ are positive constants which do
not depend on \(z,\, n\) or \(\delta\), but may depend on
\(P\left( z\right) .\)}
\end{enumerate}

\textsl{Assume that \(M(\delta) \)  satisfies the estimate}%
\begin{equation}
\int^{\frac{\pi}{2}}_{0}\log \log M(\delta )d\delta <\infty,
\label{loglog}
\end{equation}
\textsl{then the function \(P(z)\) is uniquely determined by
the formal power series \(\hat{P}(z)=\sum_{k=0}^{\infty }{p_{k}}/{z^{k+1}}\): two functions \(P_1(z)\) and \(P_2(z)\)
satisfying conditions (i) and (ii) with the same sequence
\(\{p_0,p_1,\ldots \}\)  must coincide on \( S\left( -\frac \pi
2,\frac \pi 2\right). \)}

\vspace{1.0ex}

\textbf{Corollary.} \textsl{If  the function \(M(\delta)\) in
theorem 1 is of the form}
\begin{equation}
M(\delta) = M\exp \left( \frac b{\delta ^\gamma }\right)
\label{Mdelta}
\end{equation}
\textsl{for some positive numbers \(M\),
\(b\), and \(%
\gamma \)}\textsl{ which do not depend on \(\delta\), then the
conclusion of the theorem remains true.} \vspace{2.5ex}

\textbf{Stirling's formula.} Analytic functions, satisfying
(\ref{E_S}), (\ref{k-estimates}) or (\ref {Est_KP}) have been
known in analysis for a long time.

We recall an example based upon Stirling's Formula for Euler's
Gamma Function, \(\Gamma \left( z\right) \), which can be defined
as
\[
\Gamma \left( z\right) =\int_0^\infty t^{z-1}e^{-t}dt.
\]
This function is analytic in \(S\left( -\frac \pi 2,\frac \pi
2\right),\) but it can be continued analytically to the whole of
\({\mathbb C}_z\) cut along the negative ray, that is to \(S
\left( -\pi ,\pi \right) ,\) since \(\frac 1{\Gamma \left(
z\right) }\) is an entire function with no zeros in the cut plane.

Stirling showed that
\begin{equation}
\ln \Gamma \left( z\right) =\left( z-\frac 12\right) \ln z-z+\frac
12\ln 2\pi + o\left( 1\right) ,z\rightarrow +\infty . \label{Stir}
\end{equation}
We consider the so-called \textsl{Binet function}
\begin{equation}
P(z){=}\ln \Gamma \left( z\right) -\left( z-\frac 12\right) \ln
z+z-\frac 12\ln 2\pi .  \label{Stirling1}
\end{equation}
and we will associate with this function a sequence
\(\{p_0,p_1,\ldots \}\) such that the relations (\ref{E_S}) are
valid for some choice of \(M\)  and \(a\).

According to the first Binet formula we have
\begin{equation}
P(z)=\int_0^\infty e^{-zt}F(t)dt \label{Binet1}
\end{equation}
where
\begin{equation}
F(t)=\frac{1}{t}\left(
\frac{1}{2}-\frac{1}{t}+\frac{1}{e^t-1}\right). \label{Binet2}
\end{equation}
As \(F(t)\) is analytic in the disc \(|t|<2\pi\) it can be
represented inside the disc by its Taylor series
\begin{equation}
\hat{F}(t)=\sum_{k=0}^{\infty }f_{k}t^{k}. \label{Fff}
\end{equation}
Substituting \(\hat{F}(t)\) for \(F(t)\) in (\ref{Binet1}) yields
a formal power series
\begin{equation}
\hat{P}(z)=\sum_{k=0}^{\infty }\frac{p_{k}}{z^{k+1}} ,\,\,
p_k=f_kk!.\label{Pf}
\end{equation}
Since \(F(t)\) is an even function, \(F(-t)=F(t),\) it follows
from (\ref{Pf}) that \(p_{2k-1}=0,\) for \(k=1,2,\ldots,\) and
further analysis shows that
\begin{equation}
 p_{2k-2}=\frac{B_{2k}}{2k(2k-1)},\,\,k=1,2,\ldots. \label{Bern1}
\end{equation}
Here \(B_{2k}\), \(k=0,1,\ldots,\) are the Bernoulli numbers,
which are defined as the coefficients of the Taylor series
\begin{equation}
\frac{t}{e^t-1}=1
-\frac{t}{2}+\sum_{k=2}^{\infty}B_{2k}\frac{t^{2k}}{(2k)!},
\label{Bern}
\end{equation}
 see, for example,  formulas 23.1.1-23.1.3 and 6.1.42 of
\cite{Abramowitz}. Thus \(\hat{P}(z)\) can be re-written in the
form
\begin{equation}
\hat{P}(z)=\sum_{k=1}^{\infty }\frac{B_{2k}}{2k(2k-1)z^{2k-1}},
\label{Stirling2}
\end{equation}
and this series is known as Stirling's series.\footnote{A similar
definition for Stirling's series is given in \cite{Whi}, but the
notation \(B_k\) used there refers to \((-1)^{k+1}B_{2k}\), where
\(B_{2k}\) is defined as above in (\ref{Bern}).}

Using (\ref{Binet1}) one can show that \(\hat{P}(z)\) is the
Gevrey expansion for \(P(z)\) in \(S\left( -\frac \pi 2,\frac \pi
2\right).\) Moreover, the following series of estimates is valid
\begin{align}
\left| P\left( z\right) \right| & \leq K\left( z\right) \frac{\left|
B_2\right| }{2\left| z\right| },  \nonumber \\
\left| P\left( z\right) -\sum_{k=1}^n
\frac{p_{2k-2}}{z^{2k-1}}\right| & \leq K\left( z\right)
\frac{\left| B_{2(n+1)}\right| }{(2n+2)(2n+1)\left| z\right|
^{2n+1}},
\label{Estimates_St} \\
n =1,\ldots ,\,0<|z|<+\infty &,\,z\in S\left( -\frac \pi 2,\frac
\pi 2\right),  \nonumber
\end{align}
where
\begin{equation}
K\left( z\right) ={\max\limits_{u\geq 0}}\left| \frac{z^2}{u^2+z^2}\right| ,
\label{K}
\end{equation}
(see, for example, formula 6.1.42 of \cite{Abramowitz}).

Using known asymptotics for the Bernoulli numbers, see formula
23.1.15 of \cite{Abramowitz}, we have
\begin{equation}
B_{2n+2}=\left( -1\right) ^{n+2}\frac{2\left(%
2n+2\right)!}{\left( 2\pi \right) ^{2n+2}}\left( 1+o\left(
1\right) \right) ,k\rightarrow \infty , \label{Bern2}
\end{equation}
and it follows that (\ref{Estimates_St}) is of the form given in
(\ref{E_S}) with \(a=2\pi .\) Thus, the Stirling series
(\ref{Stirling2}) is the Gevrey expansion for \(P(z)\) in the
sector \(S \left(-\frac{\pi}{2},\frac{\pi}{2} \right).\)

\vspace{1.5ex}

\textbf{Comment 2}. \textsl{The importance of Stirling's example
for this subject lies in the following observations:}

\textsl{(i) It follows from (\ref{Bern2}) that the estimate
(\ref{Estimates_St}), with the sequence }\( \{p_{0},p_{2},\ldots
\}\)\textsl{\ given by (\ref{Bern1}), ensures uniqueness in the
sector \( S\left( -\frac \pi 2,\frac \pi 2\right) , \) with
opening equal to }\(\pi .\) \textsl{This is despite the fact that
if \(z\) belongs to the boundary of \( S\left(-\frac {\pi
}{2},\frac {\pi }{2}\right)\) then the value of \(K(z)\) in the
right-hand side of (\ref {Estimates_St}) is equal to \(\infty.\)
Since the function \(K\left( z\right) \) in (\ref{Estimates_St})
can be represented as
\begin{equation}
K(z) =\begin{cases}
\hspace{6.0ex}1\hspace{6.0ex},& -\frac{\pi}{4}<\ \, \arg z\ <\frac \pi 4, \\[1.0ex]
\dfrac{1}{\sin ( 2| \arg z| ) }\ ,&{\ \ \ }\frac{\pi}{4}\leq |\arg z| < %
\frac{\pi}{ 2}.\label{1}
\end{cases}
\end{equation}
it follows that the function \(P\left( z\right) \) given by
(\ref{Stirling1}) satisfies all the conditions of Theorem 1, with
\(M(\delta) =\frac 1\delta \) and \(\sigma=0.\)}

\textsl{(ii) Substituting \({n}_{opt} (|z|)=[\pi \left| z\right|
-1]\) for \(n\) in the inequality
(\ref{Estimates_St})\footnote{This value of \(n={n}_{opt} (|z|)\)
can be guessed by minimizing the expression on the right-hand side
of (\ref{Estimates_St}) for
given \(|z|\).}  yields}%
\begin{equation}
\left| P\left( z\right) -\sum_{k=1}^{[\pi \left| z\right| -1]}\frac{p_{2k-1}%
}{z^{2k-1}}\right| \leq K\left( z\right) \frac{2\sqrt{2\pi \left| z\right| }}{%
2\pi \left| z\right| -1}e^{-2\pi \left| z\right| }. \label{Est_St}
\end{equation}
\textsl{For given fixed \(z \in S
\left(-\frac{\pi}{4},\frac{\pi}{4} \right)\)  the estimate
(\ref{Est_St}) gives the best possible accuracy when replacing
\(P\left( z\right) \) by the optimal finite sum of Stirling's
series, and it follows from (\ref{1}) that, for \(|z|>1,\) the
error is less than}
\begin{equation}
Me^{-a \left| z\right|}, \label{error}
\end{equation}
where \(M=.94891\) and \(a=2\pi. \)

\textsl{(iii) It also follows from (\ref{Binet1})-(\ref{Binet2})
 that \(\hat{P}(z)\) is a Gevrey expansion of
\(P(z)\) in the whole region \(S(-\pi,\pi)\). However, the
estimate (\ref{Estimates_St}) with \(a=2\pi\) is valid in \(
S\left( -\frac \pi 2,\frac \pi 2\right)\) only. It can be shown
that in a sector \(S\left(-{\frac {\pi }{2}}-\varepsilon,{\frac
{\pi }{2}}+\varepsilon\right)\) estimates similar to (\ref{E_S})
are valid but the constant \(a\) in these estimates is less then
\(2\pi\).  Moreover, the constant \(a\) for the sector
\(S\left(-{\frac {\pi }{2}}-\varepsilon,{\frac {\pi
}{2}}+\varepsilon\right)\) depends on \(\varepsilon .\) For this
example it can be proved that}
\begin{equation}
a(\varepsilon)=2\pi \cos\varepsilon.
\end{equation} \hfill \(\blacktriangle\)

\vspace{0.5ex}

One can derive a necessary and sufficient condition for the
uniqueness following, for example, the technique of Ahlfors'
distorsion theorem, \cite{Ahlfors}, and its complement and
refinement by Warschawski, \cite{War}. We intend to return to this
question in a later publication. Keeping in mind that we will
consider later the case in which the parameter \(a\) in
(\ref{Est_KP}) also depends on \(\delta\), and to apply our
technique to the reconstruction problem, see Watson's
representation theorem and theorem 2, pages 15-16, we prefer
another approach which we demonstrate below.

\vspace{2.0ex}

\textbf{The main Lemma.} Theorem 1 follows from the following
generalization of Lemma 2.

\vspace{0.5ex}

\textbf{Lemma 3.} \textsl{Assume that \(P\left( z\right) \) is
analytic in the right half-plane of the \(z-\)plane. Assume
further that for some
\(a>0,\) for all \(\delta \), \(%
0<\delta < \frac \pi 2\), and for all \(z\), \(z\in S\left( -\frac
\pi 2+\delta ,\frac \pi 2-\delta \right),\) the function \(P\left(
z\right) \) satisfies the following estimate}
\begin{equation}
\left| P\left( z\right) \right| <M(\delta) e^{-a\left| z\right| },
\label{Mcondition}
\end{equation}
\textsl{where \(M(\delta) , \, \log M(\delta)>1,\,\delta \in
\left(0,\frac \pi 2\right)\) satisfies condition (\ref{loglog}).
Then}
\begin{equation}
P\left( z\right) \equiv 0.
\end{equation}

\vspace{1.5ex}

\textsf{Proof of Lemma 3}. Unfortunately we cannot integrate the
expression \(e^{zt}P\left( z\right) \) along the imaginary ray as
we did in (\ref{Bo}), nor, in general, can we integrate along a
line parallel to it. This obstacle is quite typical of such
situations and we will show a way in which it may be overcome.

Given \(\theta ,\,-\frac \pi 2<\theta <\frac \pi 2,\) introduce
the Laplace transform
\begin{equation}
F_\theta \left( t\right) =\int_{l_\theta }e^{zt}P\left( z\right) dz,
\label{L-theta}
\end{equation}
where
\[
l_\theta =\left\{ z:\arg z=\theta \right\} .
\]

In what follows we will show that:

(i) The function \(F_\theta \left( t\right) \) can be continued
analytically to an entire function \(F\left( t\right) \) which
does not depend on \(\theta .\)

(ii) The function \(F\left( t\right) \) is bounded outside any
sector \( S_\delta ,\,0<\delta <\frac \pi 2,\) with angle
\(2\delta ,\) which has its apex at a point \(\frac a{2\sin \delta
},\) contains the interval \(\left( \frac a{2\sin \delta },+\infty
\right) \), and is symmetric with respect to the real line.

(iii) On the boundary \(\partial S_\delta\) of the sector
\(S_\delta \) we have \(\left| F\left( t\right) \right| \leq
\frac{2M(\delta) }a\).

(iv) The function \(F\left( t\right) \) satisfies all the
conditions of Carleman's theorem, see \cite{Ca1}, which yields
\(F\left( t\right) \equiv 0.\)

\vspace{1.5ex}

\textbf{Remark 3.} Victor Havin, to whom we have shown our result,
called our attention to the following fact. If one replaces the
condition (\ref{Mcondition}) with \(M(\delta)\) satisfying
(\ref{loglog}) by the stronger condition
\begin{equation}
\left| P\left( z\right) \right| <M\exp \left( \frac b\delta \right)
e^{-a\left| z\right| },  \label{bdelta-condition}
\end{equation}
for some \(b>0,\) then for any \(c\), \(0<c<a,\) there exists
\(h>0\) and \(M_h>0\) such that

\begin{equation}
\left| P\left( h+iy\right) \right| <M_he^{-c\left| y\right| }.
\label{c-inequality}
\end{equation}

Using the conditions of Lemma 2, the inequality
(\ref{c-inequality}) yields immediately \(P\left( z\right) \equiv
0.\)

To prove (\ref{c-inequality}) one needs only to note that if $z=x+iy$ and
\[
\frac x{\left| y\right| }=\tan \delta ,
\]
then the inequality (\ref{bdelta-condition}) can be re-written
\begin{equation}
\left| P\left( z\right) \right| <M\exp \left( \frac b{\arctan \frac x{\left|
y\right| }}\right) e^{-a\left| z\right| }.  \label{xy-condition}
\end{equation}
In turn, for
\[
\left| y\right| \gg x>0,
\]
(\ref{xy-condition}) can be rewritten in the form
\begin{equation}
\left| P\left( z\right) \right| <M^{\prime }\exp \left(
\frac{b\left| y\right| }x\right) e^{-a\left| y\right| },
\end{equation}
for some constant \(M^{\prime }>0,\) which clearly yields (\ref{c-inequality}%
), after an appropriate choice of \(x=h\), for example,
\[
h=\frac b{a-c}.
\]

Unfortunately, this beautiful argument does not cover the case of
faster growth as $\delta \rightarrow 0$, as allowed in
(\ref{loglog}), and also, for example, in (\ref{Mdelta}).\hfill
\(\blacktriangle\)

\vspace{1.5ex}
 We return now to the proof.

\vspace{1.5ex}

\textbf{Region of analyticity of \boldmath{\(F_\theta \left(
t\right) \)}}. To find the region of analyticity of \(F_\theta
\left( t\right) \) we provide the following estimate using
(\ref{Mcondition}),
\begin{equation}
\left| F_\theta \left( t\right) \right| <M(\delta) \int_{l_\theta
}\left| e^{zt}\right| e^{-a\left| z\right| }d\left| z\right| .
\label{Est-theta}
\end{equation}

Since
\begin{equation}
\left| e^{zt}\right| =e^{\left| z\right| \left( \sigma \cos \theta
-\tau \sin \theta \right) },\ \ \ t=\sigma +i\tau ,  \label{Mod}
\end{equation}
the integral (\ref{L-theta}) exists and represents an analytic function in
the half-plane
\begin{equation}
\Pi _{\theta ,a}=\left\{ t\in \mathbb{C}_t:\sigma \cos \theta
-\tau \sin \theta <a\right\} .  \label{Pi-theta}
\end{equation}
In fact, the line
\begin{equation}
L_{\theta ,a}=\left\{ t\in \mathbb{C}_t:\sigma \cos \theta -\tau
\sin \theta =a\right\}   \label{Line-a}
\end{equation}
divides the \(t-\)plane into two half-planes, and \(\Pi _{\theta
,a}\) is the half-plane containing the origin \(t=0,\) see Figure
1, where we assume that \(0<\theta <\frac \pi 2\).

\begin{wrapfigure}[12]{l}{0.5\linewidth}
\begin{minipage}[H]{1.0\linewidth}

\vspace*{-1.5ex}
\begin{minipage}[H]{1.0\linewidth} \epsfig{file=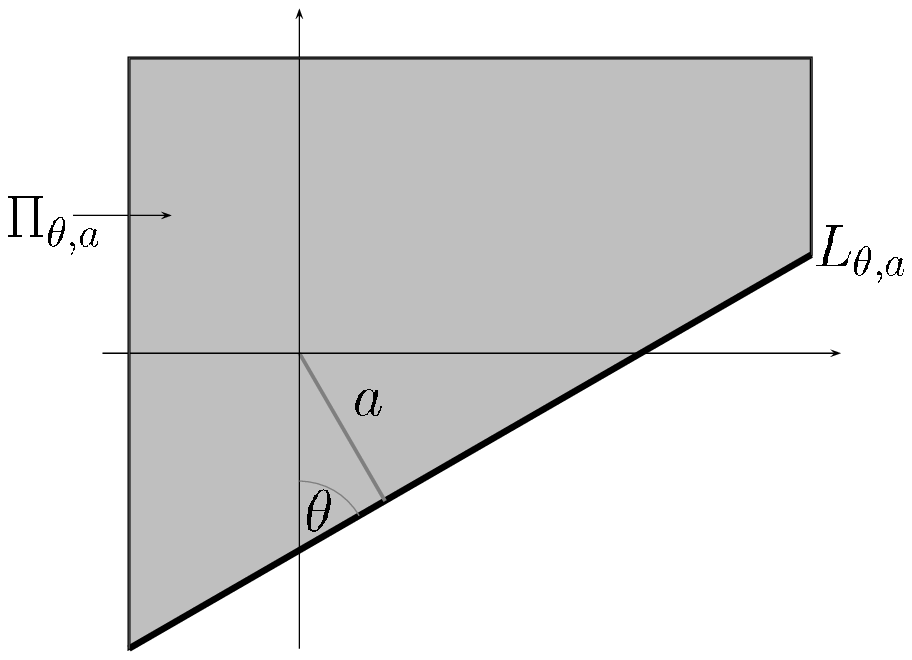,width=0.95\linewidth,clip=}
%\caption{The line $L_{\theta ,a}$ and the
%half-plane $\Pi _{\theta ,a}$}
\label{F1}%
\end{minipage}

 \vspace*{+1.5ex}\hspace*{1.5ex}
\begin{minipage}[H]{1.0\linewidth}
\textbf{Figure 1}:\ \ The line \(L_{\theta ,a}\) and \\ the
half-plane \(\Pi _{\theta ,a}\)
\end{minipage}
\end{minipage}
\end{wrapfigure}
\vspace{0.5ex}

Moreover, from (\ref{Est-theta}) and (\ref{Mod}) we have for
\(t\in \Pi _{\theta ,a}\) the estimate
\begin{equation}
\left| F_\theta \left( t\right) \right| <M(\delta) \frac
1{a-\left( \sigma \cos \theta -\tau \sin \theta \right)
}\label{F-est1}
\end{equation}
where
\[
\delta =\left\{
\begin{array}{c}
\frac \pi 2-\theta ,\theta >0 \\
\frac \pi 2+\theta ,\theta <0
\end{array}
\right \}
\]

In what follows we will show that the functions \(F_\theta \left(
t\right),\,\left|\theta \right| <\frac \pi 2,\) are all elements
of a single analytic function \(F\left( t\right) \) which does not
depend on \(\theta .\)

\vspace{1.5ex}

\textbf{Relationship between \boldmath{\(F_\theta \left(
t\right)\)}\( \) and \boldmath{\(F_{-\theta} \left( t\right)\)}\(
\).} Given \(\delta ,\,0<\delta <\frac \pi 2\), consider the pair
of functions \(F_{\frac \pi 2-\delta }\left( t\right) \) and
\(F_{-\frac \pi 2+\delta }\left( t\right) \) of the form (\ref
{L-theta}) which are analytic in the half-planes \(\Pi _{\frac \pi
2-\delta ,a}\) and \(\Pi _{-\frac \pi 2+\delta ,a}\) respectively.
Now
\begin{equation}
\Pi _{\frac \pi 2-\delta ,a}=\left\{ t\in \mathbb{C}_t:\sigma \sin
\delta -\tau \cos \delta <a\right\}   \label{Pi-1}
\end{equation}
and
\begin{equation}
\Pi _{-\frac \pi 2+\delta ,a}=\left\{ t\in \mathbb{C}_t:\sigma
\sin \delta +\tau \cos \delta <a\right\}   \label{Pi-2}
\end{equation}
and these half-planes are symmetric with respect to the real line
of the \(t-\)plane.

\begin{wrapfigure}[13]{l}{0.6\linewidth}
\begin{minipage}[H]{0.9\linewidth}
\vspace*{1.8ex} \hspace*{2.0ex}
 \epsfig{file=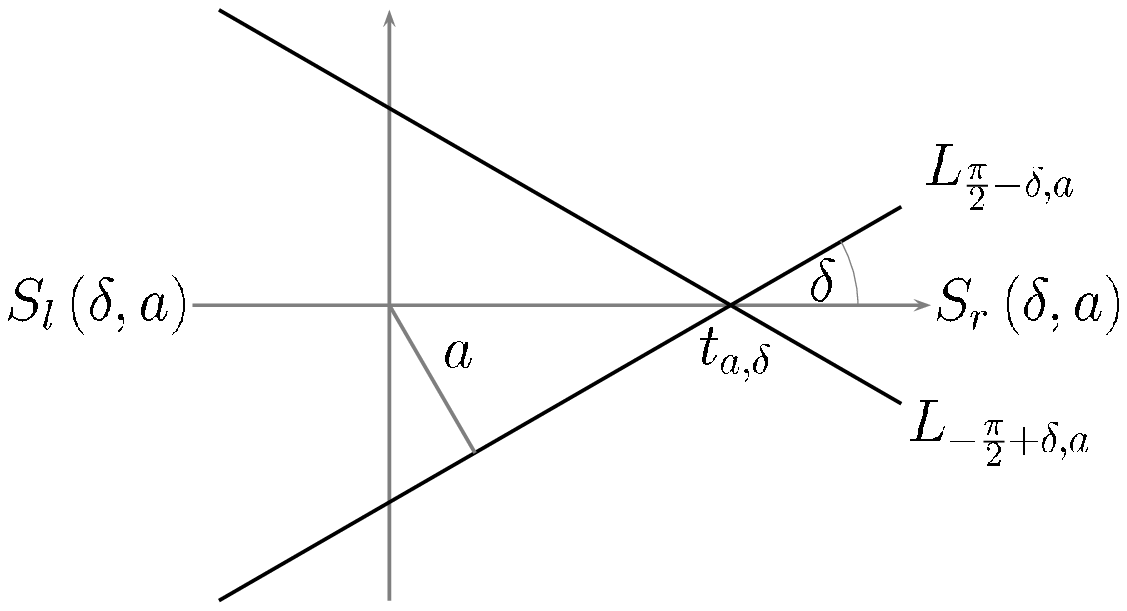,width=1.0\linewidth,clip=}
%\caption{The lines {$L_{\mp\frac \pi 2\pm\delta,a}$} and the
%sectors {$S_l\left( \delta ,a\right) $} and {$S_r \left( \delta
%,a\right) $} }
\label{F2}%
\end{minipage}\\
\hspace*{3.0ex}
\begin{minipage}[H]{0.8\linewidth}%
\textbf{Figure 2}:\ \ The lines \(L_{\mp\frac \pi
2\pm\delta,a}\)\\
and the sectors \(S_l( \delta ,a)\) and \(S_r ( \delta ,a) \).
\end{minipage}
\end{wrapfigure}
On Figure 2 we show the two lines \(L_{\frac \pi 2-\delta ,a}\)
and \( L_{-\frac \pi 2+\delta ,a}\) crossing the positive ray of
the \(t-\)plane at the point
\begin{equation}
t_{\delta ,a}=\frac a{\sin \delta },  \label{apex}
\end{equation}
at angles \(\delta \) and \(-\delta \) respectively.

These lines define two sectors, \(S_l\left( \delta ,a\right) \)
(the left-hand sector) and \(S_r\left( \delta ,a\right) \) (the
right-hand sector), both with their apexes at \( t_{\delta ,a}\)
and with angle equal to \(2\delta \).

\hspace{3.0ex}

Set

\begin{equation*}
S_1\left( \delta ,a\right) =\Pi _{\frac \pi 2-\delta ,a}\cap \Pi
_{-\frac \pi 2+\delta ,a}
\end{equation*}
and
\begin{equation}
S_2\left( \delta ,a\right) =\Pi _{\frac \pi 2-\delta ,a}\cup \Pi _{-\frac
\pi 2+\delta ,a}.  \label{Pi-union}
\end{equation}

The left-hand sector \(S_l\left( \delta ,a\right) \) is exactly
\(S_1\left( \delta ,a\right) ,\) while the closure of the
right-hand sector \(S_r\left( \delta ,a\right) \) is the
complement of \(S_2\left( \delta ,a\right) .\) (See Figures 3 and
4.)

\begin{minipage}[H]{1.0\linewidth}
\begin{minipage}[H]{1.0\linewidth}
\begin{minipage}[H]{0.44\linewidth}
 \epsfig{file=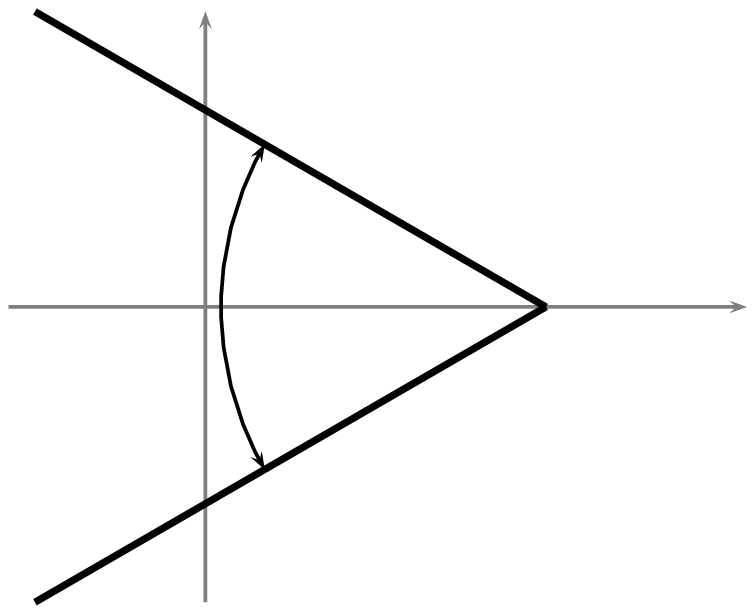,width=0.78\linewidth,clip=}
\end{minipage}\hspace*{1.0ex}%
\begin{minipage}[H]{0.48\linewidth}
 \epsfig{file=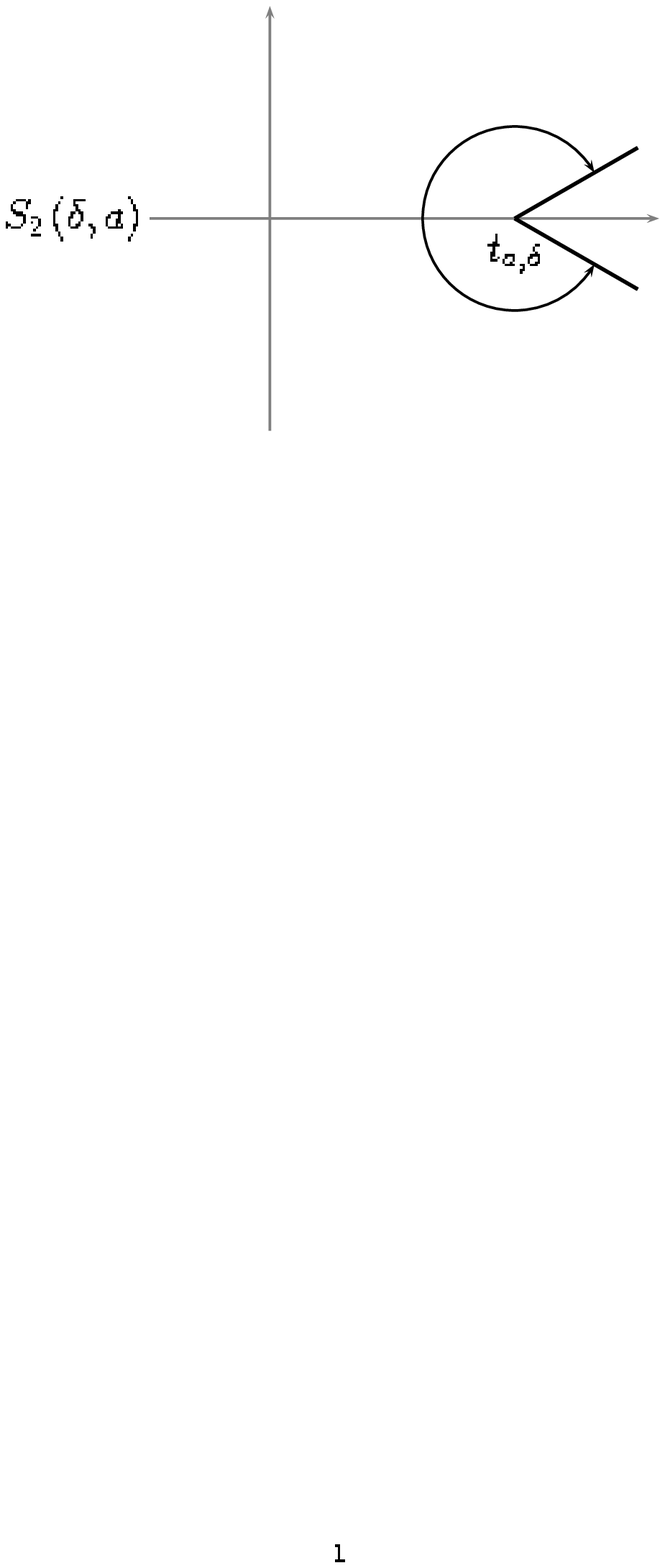,width=0.99\linewidth,clip=}
\end{minipage}
\end{minipage}

\vspace*{2.0ex}%
{} \hspace*{-3.0ex} \textbf{Figure 3}:\ \ The sector
\(S_{1}(\delta,a)\).\hspace{5.0ex} \textbf{Figure 4}:\ \ The
sector \(S_{2}(\delta,a)\).
\end{minipage}

\vspace{5.0ex}

\noindent
Using the representation (%
\ref{L-theta}) and Cauchy's theorem yields
\[
F_{\frac \pi 2-\delta }\left( t\right) =F_{-\frac \pi 2+\delta
}\left( t\right),\,t\in S_1\left( \delta ,a\right) .
\]
Thus, \(F_{\frac \pi 2-\delta }\left( t\right) \) and \(F_{-\frac
\pi 2+\delta }\left( t\right) \) can be considered as elements of
the function \(F\left( t,\delta \right) \) where
\begin{equation}
F( t,\delta) =
\begin{cases}
F_{\frac \pi 2-\delta }\left( t\right) ,&t\in \Pi _{\frac \pi
2-\delta ,a}
\\[1.5ex]
F_{-\frac \pi 2+\delta }\left( t\right) ,& t\in \Pi _{-\frac \pi
2+\delta ,a}
\end{cases}
 .  \label{F-delta}
\end{equation}
and this function is analytic in \(S_2\left( \delta ,a\right) .\)

However, \(S_2\left( \delta ,a\right) \) still does not give the
maximal region of analyticity of the function \(F\left( t\right)
\) referred to in (i) of page 10.

\vspace{1.5ex}

\textbf{The maximal region of analyticity}. We require the following lemma.

\vspace{1.5ex}

\textbf{Lemma 4}. \textsl{The functions }\(F\left( t,\delta \right) \)\textsl{%
\ given by (\ref{F-delta}) are all elements of a single analytic
function }\( F\left( t\right) \)\textsl{\ which does not depend on
}\(\delta \)\textsl{\ and which is an entire function. }

\vspace{1.5ex}

\textsf{Proof}. Consider \(\delta ^{\prime }>\delta ^{\prime
\prime }>0\) and the corresponding functions \(F\left( t,\delta
^{\prime }\right) \) and \(F\left( t,\delta ^{\prime \prime
}\right) .\)

We have
\[
S_r\left( \delta ^{\prime },a\right) \supset S_r\left( \delta ^{\prime
\prime },a\right) ,
\]
and so
\[
{\mathbb{C}_t}\backslash S_r\left( \delta ^{\prime },a\right) \subset {\mathbb{C}_t}%
\backslash S_r\left( \delta ^{\prime \prime },a\right) .
\]

Using the same argument as that used above we see that
\[
F\left( t,\delta ^{\prime }\right) =F\left( t,\delta ^{\prime
\prime }\right),\,t\in {\mathbb{C}_t}\backslash S_r\left( \delta
^{\prime },a\right) ,
\]
and \(F\left( t,\delta ^{\prime \prime }\right) \) can be
considered as an analytical continuation of \(F\left( t,\delta
^{\prime }\right) \) to the larger region
\({\mathbb{C}_t}\backslash S_r\left( \delta ^{\prime \prime
},a\right) \).

The expression (\ref{apex}) for \(t_{\delta ,a}\) shows that
\begin{equation}
 t_{\delta ,a}\rightarrow \infty ,\, \delta \rightarrow 0,
 \label{infinity}
\end{equation}
which yields
\[
\bigcap _{{0<\delta <\frac \pi 2}}S_r\left( \delta ^{\prime
},a\right) =\varnothing .
\]

Thus the functions \(F\left( t,\delta \right) ,\,0<\delta <{\frac
\pi 2,}\) are all elements of a single function \(F\left( t\right)
,\) which is an entire function, and so statement (i) is verified.
\hspace*{0.01ex}\hfill \(\blacktriangle \)

\vspace{1.5ex}

\textbf{Estimates for \boldmath{\(F\left( t\right)\)}\( \)
 in the \boldmath{\(t-\)%
plane.}} Using the estimate (\ref{F-est1}) for the region
\begin{equation}
\bar{\Pi}_{\theta ,\frac a2}=\Pi _{\theta ,\frac a2}\cup L_{\theta ,\frac
a2},  \label{Pi-L}
\end{equation}
where \(\Pi _{\theta ,\frac a2}\ \)and \(L_{\theta ,\frac a2}\) are defined by (%
\ref{Pi-theta}) and (\ref{Line-a}) respectively, with \(a\)
replaced by \(\frac a2\), yields
\begin{equation}
\left| F_\theta \left( t\right) \right| <\frac{2M(\delta) }a,\,t\in \bar{\Pi}%
_{\theta ,\frac a2}.  \label{F-est2}
\end{equation}

Hence, from (\ref{F-delta}) we have
\begin{equation}
\left| F\left( t,\delta \right) \right| <\frac{2M(\delta)
}a,\,t\in S_2\left( \delta ,\frac a2\right) . \label{F-est3}
\end{equation}

We now apply this estimate for the function \(F\left( t\right)\)
given in (i).

Now in the region \(S_2\left( \delta ,\frac a2\right) \) the
function \(F\left( t,\delta \right) ,\) and hence also the
function \(F\left( t\right) \), is bounded, which verifies (ii).
Moreover, it follows from (\ref{F-est3}) that \(F\left( t\right)
\) satisfies the estimate
\begin{equation}
\left| F\left( t\right) \right| \leq\frac{2M(\delta) }a.
\label{F-est4}
\end{equation}
in the sector \(S_2\left( \delta ,\frac a2\right) \), and hence
also on its boundary. This proves (iii).

The inequality (\ref{F-est4}) implies
\begin{equation}
|F(|t|e^{\pm i\delta})|<\frac{2M(\delta)}{a}\label{F-est5}
\end{equation}
which shows that the conditions of the following theorem of
Carleman are satisfied.

\vspace{1.5ex}

\textbf{Carleman's Theorem}, see \cite{Ca1}. {\it Let
\(M(\varphi)\) be a positive function (finite or infinite), for
which \(\log M(\varphi)>1\) and the following integral,
\begin{equation}
\int_0^{2\pi}\log\log M(\varphi)d\varphi,
\end{equation}
exists. Every entire function \(f(z)\) which satisfies the
inequality
\begin{equation}
|f(z)|<M(\varphi), \, \varphi =\arg z,\,\;0<\varphi <2\pi,
\end{equation}
is a constant.\footnote{If the function \(M(\varphi)\) satisfies
the stronger condition (\ref{Mdelta}) this fact follows
immediately from the Phragm\'{e}n-Lindel\"{o}f theorem.}}

\vspace{1.5ex}

Thus \(F(t)=C\) for some constant \(C.\) To derive (iv) it remains
to note that since for every \( \delta \) we have
\[
F\left( t,\delta \right) \rightarrow 0,\,t\rightarrow -\infty ,
\]
it follows that
\[
F\left( t\right) \rightarrow 0,\,t\rightarrow -\infty .
\]
Thus
\[
F\left( t\right) \equiv 0.
\]
This proves (iv) and then
\[
P\left( z\right) \equiv 0
\]
follows immediately.\hfill \(\blacktriangle \)

\vspace{1.5ex}

\textbf{Comment 4.} (i) Our proof shows that it is possible to
extend our result to the case in which \(a\) depends on \(\delta
.\) From (\ref{apex})
 it follows that (\ref{infinity}) holds only
if
\begin{equation}
\frac{a(\delta )}{\delta}\rightarrow +\infty \text{ as } \delta
\rightarrow 0. \label{a-condition}
\end{equation}
Thus (\ref{a-condition}) is a necessary condition for the
uniqueness result to hold.

(ii) Carleman's result was further developed by N. Levinson and
N. Sjoberg (see the survey by M. Sodin, \cite{Sod}, for details)
who obtained a more general result, from which there follows, in
particular, an extension of Carleman's theorem in which the sector
is replaced by a half-strip. We believe that this extension may be
of assistance in further investigating the restoration problem
which we discuss below.

\vspace{1.5ex}

 Watson's uniqueness theorem shows that there is
the possibility of reconstructing, in principle, a function
\(P\left( z\right) \) satisfying the conditions (i) and (ii) of
this theorem from the formal power series \(\hat{P}\left( z\right)
\). This possibility can be realized using the Borel summation
method. The corresponding result is known as Watson's
representation Theorem. The proof of this Theorem can be extracted
from \cite{Wat1}, Section 9, though this statement has not been
represented there in the form of a theorem. See also \cite{Wat2},
page 68, formulae (20), (21) and (22), and \cite{Hardy}, section
8.11.

\vspace{1.5ex}

\textbf{Watson's representation (reconstruction) theorem}.
\textsl{Assume that the function }$P\left( z\right) $
\textsl{satisfies all the conditions of Watson's uniqueness
theorem in the sector}
\begin{equation}
S\left( -\frac \pi 2-\varepsilon ,\frac \pi 2+\varepsilon \right),
\label{region1}
\end{equation}
\textsl{except \(|z|<\sigma,\) for some \(\varepsilon,\,0<\varepsilon <\frac \pi 2.\)}

\textsl{Set}

\begin{equation}
\hat{F}\left( t\right) =\sum_{n=0}^\infty \frac{p_n}{n!}t^n.  \label{F}
\end{equation}

\textsl{Then}
\begin{enumerate}
\item[\textup{(i)}] \textsl{the formal power series} $\hat{F}\left(
t\right) $ \textsl{is} \textsl{the Taylor series for a function
}$F\left( t\right) ,$\textsl{\ which is analytic in a disc \(D_a\)
of radius }$a$\textsl{, centered at the origin.}
\item[\textup{(ii)}]  \textsl{the
function \(F\left( t\right)\) can be continued analytically from
the disc \(D_a\) to the region}
\begin{equation}
D_a \cup \left\{ t\in \mathbb{C}%
_t:\left| \arg t\right| <\varepsilon \right\} .  \label{Region1}
\end{equation}

\item[\textup{(iii)}]\textsl{ the function }$P\left( z\right) $\textsl{\ can be
represented in the form}
\begin{equation}
P\left( z\right) =\int_0^{+\infty }F\left( t\right)
e^{-zt}dt,\,\,z \in S\left( -\frac \pi 2,\frac \pi 2 \right), \label{LF}
\end{equation}
\end{enumerate}
\textsl{and the integral is absolutely convergent.}
\vspace{1.5ex}

\textbf{Remark 5}. Two problems immediately suggests themselves:
(a) what happens to the region (\ref{Region1}) of analyticity of 
\(F(t)\) as \(\varepsilon \rightarrow 0\)? (b) given that \(P(z)\) satisfies
(\ref{E_S}) in \(S\left( -\frac \pi 2-\varepsilon ,\frac \pi 2+\varepsilon \right)\) with fixed 
\(a\) and \(\varepsilon\), can the region of analyticity be improved?  

The problem (a) was answered by F. Nevanlinna in 1918, 
see \cite{Nev}, even for the more general case of Gevrey expansions 
of an arbitrary order \(k\). The improvement of Watson's representation
theorem (problem (b)) was obtained 62 years later \footnote{Ramis' book,
\cite{Ramis}, explains why Watson's theory was further developed
only recently. } by Alan D. Sokal who applied it 
to the perturbation expansion in the \({\phi
_2}^4\) quantum field theory, see \cite{Sok}. 
Sokal's improvement gives a necessary and sufficient characterization for a large class
of Borel-summable functions.   
It seems that one of the first re-appearances of
Watson's theorems was in connection with the theory of anharmonic
oscillators, see \cite{GGS}.\hfill
\(\blacktriangle\)

\vspace{1.5ex}

Our theorem 1 suggests a corresponding result for a function
\(P\left( z\right) \), satisfying the inequalities (\ref{Est_KP}) in
the critical sector \(S\left( -\frac \pi 2,\frac \pi 2\right) \)
with \(M(\delta) \)
satisfying (\ref{loglog}). Consider the simplest case, for which \(%
M(\delta) \) is a positive constant. For this case we have the
extension of Watson's uniqueness theorem alluded to in Lemma 2.
This suggests the following representation theorem for \(P\left(
z\right) .\)

Let \(D_a\) be the disc of radius \(a\) centered at \(t=0\) and let
\(L_a^{+}\) be the half strip of the form
\begin{equation}
L_a^{+}=\left\{ z:\Re z>0,\left| \Im z\right| <a\right\} .
\label{Strip}
\end{equation}

\vspace{1.5ex}

\textbf{Theorem (F. Nevanlinna)}. \textsl{Assume that the function} \(P\left(
z\right) \) \textsl{is analytic in the right half-plane of the
}\(z-\)\textsl{plane and satisfies the conditions (i) and (ii) of
Theorem 1. Assume further that the function} \(M(\delta)\)
\textsl{in (\ref{Est_KP}) is a bounded function in the interval
\(\left(0,\frac \pi 2 \right)\)}. \textsl{\ Let \(F\left( t\right)
\) be a function analytic in the disc \(D_a\) which is
represented in \(D_a\)
as a sum of the power series in \(t\) given by (\ref{F}). Then \(F\left( t\right) \)%
can be analytically continued from the disc \(D_a\) to the
``half-strip''}
\begin{equation}
\mathcal{D}_a=D_a\cup L_a^{+}  \label{Region2}
\end{equation}
\textsl{with the following bound in \(\mathcal{D}_a\)}
\begin{equation}
|F(t)| \leq {K_{a'}}\exp{(\sigma |t|)}
\end{equation}
\textsl{in every subregion \(\mathcal{D}_a',\,\, a'<a\), where the positive constant \(K_{a'}\)
depends on \(a'\).
Moreover, for \(\Re z>\sigma\) the
function \(P\left( z\right) \)
 can be represented in the
form (\ref{LF}) as the Laplace transform of \(F\left( t\right) \)
in which the integral is absolutely convergent.}

\vspace{1.5ex}

It is our aim to extend this result to more general uniqueness
classes of functions satisfying (%
\ref{Est_KP}) with \(M(\delta) \) unbounded. We conjecture that it
is possible to restore \(P\left( z\right) \ \) satisfying the
conditions of Theorem 1 from its Gevrey expansion $\hat{P}\left(
z\right) $ using Borel summation as above if (\ref{loglog}) is
replaced by the stronger condition
\[
M(\delta) <M\exp \left( \frac b\delta \right) .
\]
For classes of functions satisfying (\ref{loglog}) or
(\ref{Mdelta}) we need to look for an alternative summation
method. The part (iii) of our comment 4 and the recent publication
\cite{BLS} may be of assistance here.

For a discussion and a further extension of the results 
of F. Nevanlinna and A. D. Sokal, see
\cite{GG1}.

\vspace{1.5ex}

\textbf{Conclusion.} Two points have motivated this work. The
first is the discovery by Ramis and Sibuya that, roughly speaking,
if a formal power series satisfies an analytical (non-linear)
differential equation then there exists \(k>0\) such that for each
sector with (critical) opening \(\frac \pi k\) there exists a
regular solution of the equation for which the formal solution is
a Gevrey expansion of order k. Moreover, they proved the
uniqueness of such a solution, and it is this which lead to a
re-examination of Watson's theorem for critical openings of the
sector. If P(z) is the regular solution of the Ramis-Sibuya
theorem then we believe that \(M(\delta)\) will satisfy (12) and,
probably, the stronger condition (13). However, this is an open
question.

 The second point relates to the calculation of the best estimate
of a function (and of the associated error) using an optimal
finite sum of its Gevrey expansion. This can be done using
estimates for \(M(\delta)\) and \(a^{-1}(\delta)\) in a given
sector just as was done above in (\ref{error}) for Stirling's
series.

\vspace{1.5ex}

\textbf{Acknowledgment}. The authors wish to express their
gratitude to Nick Garnham and to Peter Jones for their
encouragement and support. The authors are grateful to Victor
Katsnelson and David Lucy for numerous helpful discussions and
also to Iossif Ostrovski who read very carefully the text and
whose corrections and suggestions allowed us to improve the paper.

\noindent David W.H. GILLAM, \(<\)dgillam@swin.edu.au\(>\)
\vspace{1.0ex}
\newline
 \noindent Vladimir P. GURARII,
\(<\)vgurarii@swin.edu.au\(>\)
\vspace{1.5ex}
\newline
\noindent Faculty of Engineering and Industrial Sciences
\newline Swinburne University of Technology, \newline PO Box 218
Hawthorn VIC 3122 Melbourne, Australia

\end{document}